\documentclass{amsproc}

\usepackage{ifthen}
\newboolean{article}
    \setboolean{article}{true}
\newboolean{report}
\newboolean{book}
\newboolean{letter}
\newboolean{german}
\newboolean{nobaselinestretch}
    \setboolean{nobaselinestretch}{true}
\newboolean{nosectionappendix}
    \setboolean{nosectionappendix}{true}
\newboolean{oldtoc}
\newboolean{nosectionequn}
\newboolean{notheorem}
\usepackage{theorem}

\ifthenelse{\boolean{german}}{
    \usepackage{german}}{}

\usepackage[latin1]{inputenc}

\usepackage{amsmath}
\usepackage{amssymb}
\usepackage[mathscr]{eucal}

\ifthenelse{\boolean{nobaselinestretch}}{}{
    }

\newenvironment{env}[2]{\begin{#1}#2\end{#1}}{}
    \newcommand{\beq}[1]{\begin{env}{equation}{#1}}
    \newcommand{\beqn}[1]{\begin{env}{equation*}{#1}}
    \newcommand{\bal}[1]{\begin{env}{align}{#1}}
    \newcommand{\baln}[1]{\begin{env}{align*}{#1}}
    \newcommand{\bga}[1]{\begin{env}{gather}{#1}}
    \newcommand{\bgan}[1]{\begin{env}{gather*}{#1}}
    \newcommand{\bflal}[1]{\begin{env}{flalign}{#1}}
    \newcommand{\bflaln}[1]{\begin{env}{flalign*}{#1}}
    \newcommand{\bmu}[1]{\begin{env}{multline}{#1}}
    \newcommand{\bmun}[1]{\begin{env}{multline*}{#1}}
    \newcommand{\bsp}[1]{\begin{env}{split}{#1}}

    \newcommand{\eeq}{\end{env}}
    \newcommand{\eeqn}{\end{env}}
    \newcommand{\eal}{\end{env}}
    \newcommand{\ealn}{\end{env}}
    \newcommand{\ega}{\end{env}}
    \newcommand{\egan}{\end{env}}
    \newcommand{\eflal}{\end{env}}
    \newcommand{\eflaln}{\end{env}}
    \newcommand{\emu}{\end{env}}
    \newcommand{\emun}{\end{env}}
    \newcommand{\esp}{\end{env}}

\newcommand{\lf}{\vspace{2ex}}

\renewcommand{\bf}[1]{\textbf{#1}}
\renewcommand{\it}[1]{\textit{#1}}

\renewcommand{\sf}[1]{\textsf{#1}}

\renewcommand{\tt}[1]{\texttt{#1}}
\newcommand{\hl}[1]{\bf{\it{#1}}}

\newcommand{\mbf}[1]{\mathbf{#1}}

\newcommand{\cmc}[1]{\mathcal{#1}}
\newcommand{\eus}[1]{\mathscr{#1}}
\newcommand{\euf}[1]{\mathfrak{#1}}
\newcommand{\bb}[1]{\mathbb{#1}}
\newcommand{\nbd}[1]{$#1$\nobreakdash--}
\newcommand{\ol}[1]{\overline{#1}}

\newcommand{\vt}{\vartheta}

\newcommand{\vp}{\varphi}
\newcommand{\om}{\omega}
\newcommand{\Om}{\Omega}

\newcommand{\bfam}[1]{\bigl(#1\bigr)}
\newcommand{\Bfam}[1]{\Bigl(#1\Bigr)}
\newcommand{\AB}[1]{\langle#1\rangle}
\newcommand{\bAB}[1]{\bigl\langle#1\bigr\rangle}

\newcommand{\CB}[1]{\{#1\}}
\newcommand{\bCB}[1]{\bigl\{#1\bigr\}}
\newcommand{\BCB}[1]{\Bigl\{#1\Bigr\}}
\newcommand{\SB}[1]{[#1]}
\newcommand{\bSB}[1]{\bigl[#1\bigr]}

\newcommand{\set}[2][]{
    \ifthenelse{\equal{#1}{}}{
        \CB{#2}}{
        \CB{#1~|~#2}}}
\newcommand{\bset}[2][]{
    \ifthenelse{\equal{#1}{}}{
        \bCB{#2}}{
        \bCB{#1~|~#2}}}
\newcommand{\Bset}[2][]{
    \ifthenelse{\equal{#1}{}}{
        \BCB{#2}}{
        \BCB{#1~\big|~#2}}}

\DeclareMathOperator{\id}{\normalfont\sf{id}}

\DeclareMathOperator{\alg}{\normalfont\sf{alg}}

\newcommand{\C}{\bb{C}}

\newcommand{\E}{\bb{E}}

\newcommand{\R}{\bb{R}}

\newcommand{\cA}{\cmc{A}}
\newcommand{\cB}{\cmc{B}}

\newcommand{\sB}{\eus{B}}

\newcommand{\sN}{\eus{N}}

\newcommand{\ei}{\euf{i}}

\newcommand{\ep}{\euf{p}}

\newcommand{\eF}{\euf{F}}

\newcommand{\eH}{\euf{H}}

\newcommand{\U}{\mbf{1}}

\newcommand{\G}{\Gamma}

\ifthenelse{\boolean{nosectionequn}}{}{
    \numberwithin{equation}{section}
    }

\ifthenelse{\boolean{article}\or\boolean{letter}\or\boolean{nosectionequn}}{
    \setboolean{nosectionappendix}{true}}{}
\ifthenelse{\boolean{nosectionappendix}}{}{
    \renewcommand{\appendix}{
        \chapter*{\appendixname}
        \addcontentsline{toc}{chapter}{\appendixname}
        \renewcommand{\thesection}{\Alph{section}}
        \setcounter{section}{0}}}
   
\ifthenelse{\boolean{report}\or\boolean{book}}{
    }{}

\newcounter{s}[section]
\newcommand{\stru}[1]{\lf\noindent\stepcounter{s}\bf{\thesection.\arabic{s}~#1}}

\usepackage[matrix,arrow]{xy}

\begin{document}


\bibliographystyle{alpha}
\title{
Independence and Product Systems}

\author{Michael Skeide}

\date{January 2003}



\begin{abstract}
Starting from elementary considerations about independence and Markov processes in classical probability we arrive at the new concept of conditional monotone independence (or operator-valued monotone independence). With the help of product systems of Hilbert modules we show that monotone conditional independence arises naturally in dilation theory.
\end{abstract}

\maketitle

\vfill

\noindent
One of the most fundamental concepts in probability theory --- in the sequel, we say \hl{classical probability} --- is \hl{independence}. The most fundamental concept dealing with non-independent random variables is the \hl{Markov property}.  In order to underline the formal aspect of the Markov property being a generalization of independence (replacing expectations by conditional expectations) we will say \hl{conditional independence} (which operator algebraists would prefer to call operator-valued independence).

The expectation is a normalized positive linear functional (i.e.\ a \hl{state}) on the \hl{commutative} algebra of random variables. \hl{Noncommutative} or \hl{quantum probability} is a noncommutative generalization of classical probability designed to include also quantum physical applications. A \hl{quantum probability space} is, therefore, a pair $(\cA,\vp)$ of a unital \nbd{*}algebra (more specifically, a (pre--)\nbd{C^*}algebra or a von Neumann algebra) $\cA$ with a \hl{normalized} (i.e.\ $\vp(\U)=1$) \hl{positive} (i.e.\ $\vp(a^*a)\ge0$) linear functional $\vp\colon\cA\rightarrow\C$. In this context a classical probability space $(\Om,\eF,P)$ corresponds to the commutative quantum probability space $(L^\infty(\Om),\vp=\int\bullet\,dP)$.

Classical independence is \hl{symmetric}, i.e.\ if $X_1$ is independent of $X_2$ (conditioning $X_1$ on $X_2$ does not change probabilities of $X_1$), then $X_2$ is independent of $X_1$, too. In Section \ref{indsec} we discuss an example where this symmetry is not desirable. However, within the category of commutative quantum probability spaces there is only one independence, namely, that one which gives back classical independence. Therefore, if we want to model the mentioned example, then we are forced to leave the category of commutative quantum probability spaces even for classical random variables. The only noncommutative independence which earns the name independence and is not symmetric is \hl{monotone independence}; see Section \ref{indsec}.

In Section \ref{condindsec} we discuss the case of conditional independence, which is the basis for our discussion of Markov processes, later on. As a new concept we introduce \hl{conditional monotone independence} (or operator-valued monotone independence). As long as we are dealing with classical random variables, both conditional independence and conditional monotone independence are possible.  However, if we start with noncommutative random variables from the beginning, then the constructions for conditional indpendence are no longer possible, in general. Only in few exceptional cases (for instance, when the involved quantum probability spaces are algebras isomorphic to some algebra $\sB(H)$ of all bounded operators on some Hilbert space $H$) quantum conditional independence still exists.

In Section \ref{dilsec}, finally,  we point out with the help of \hl{product systems of Hilbert modules} that conditional monotone independence arises in every \hl{reversible quantum dynamical system} which is a \hl{dilation} of an \hl{irreversible quantum dynamical system}. (Conditional) monotone independence is, therefore, one of the most natural quantum independences possible.

\section{Independence}\label{indsec}

Let $X_1,X_2$ stand for two observable quantities. These can be \hl{classical random variables} (i.e.\ real-valued measurable functions on some probability space) or \hl{quantum random variables} (i.e.\ self-adjoint operators on some Hilbert space with a distinguished vector state $\vp=\AB{\xi.\bullet\xi}$).

To avoid technical discussions, we assume that the $X_i$ are bounded. (This is not a serious restriction, because on the classical side everything is determined when known on indicator functions, and on the quantum side by the spectral theorem everything is known when known on projections.)

Let us start with the case of two classical random variables. Everything one can say in the classical description can be formulated in terms of (marginal and joint) distributions.
$$
\begin{array}[t]{ll}\vspace{.5ex}
X_i~\rightsquigarrow~L^\infty(\R,\mu_i)~\subset~L^1(\R,\mu_i)
&\text{{marginal distribution}}
\\\vspace{1ex}
\E\SB{f(X_i)}=\int f(x)\,\mu_i(dx)
&
\\\vspace{.5ex}
X_1,X_2~\rightsquigarrow~L^\infty(\R\times\R,\mu_{12})
&\text{{joint distribution}}
\\
\E\SB{f(X_1,X_2)}=\int f(x_1,x_2)\,\mu_{12}(dx_1,dx_2)
&
\end{array}
$$
The quantum description of the classical case is the following reinterpretation.
$$
\begin{array}[t]{ll}\vspace{.5ex}
\cA_i~=~L^\infty(\R,\mu_i)~\subset~\sB(L^2(\R,\mu_i))
&\text{{algebras}}
\\\vspace{1ex}
\vp_i(f)=\int f(x)\,\mu_i(dx)=\AB{1,f1}
&\text{with a {state}}
\\\vspace{.5ex}
\cA_{12}~=~L^\infty(\R\times\R,\mu_{12})
&\text{{a ``product''} of algebras}
\\
\vp_{12}(f)=\int f(x_1,x_2)\,\mu_{12}(dx_1,dx_2)=\AB{1,f1}
&\text{with a state on the``product''}
\end{array}
$$
$X_1$ and $X_2$ are \hl{independent}, if $\mu_{12}=\mu_1\otimes\mu_2$. Then
\baln{
\cA_{12}~~~~~&=~~~~~~~~~~\cA_1\otimes\cA_2&\vp_{12}~&=~\vp_1\otimes\vp_2
\\
L^2(\R\times\R,\mu_{12})~&=~L^2(\R,\mu_1)\otimes L^2(\R,\mu_2)&1~&=~~~1\otimes 1.
}\ealn

There are many possible interpretations of independence. We mention only two. The \hl{intrinsic} interpretation: There are two random variables and one does not influence the other. In probabilistic language, conditioning on one of them does not influence knowledge of the other. In physical language, measuring one of them does not change the state of the other. We say intrinsic, because we start with a given system (probability space) and we want to decide, whether two random variables of that system are independent or not. The \hl{constructive} interpretation: Given two (quantum or classical) probability spaces with one random variable each, we construct a joint probability space on which both random variables can be described simultaneously in such a way that they influence each other as little as possible. In the preceding discussion, independence, formulated intrinsically as factorization of probabilities, gives rise to a tensor product. The tensor product, in turn, can also serve a prescription how to construct a joint probability space from the marginal distributions.

Let us start with the constructive interpretation. Schürmann \cite{MSchue95b} shows that the \hl{tensor product} of unital \hl{commutative} {quantum} {probability} {spaces} $(\cA_i,\vp_i)$ is the only prescription, how to compute mixed moments $\vp_{12}(X_1^nX_2^m)$, such that certain axioms are fulfilled. These axioms are properties derived from classical independence like associativity of the construction and the requirment that functions of independent random variables should again be independent (\hl{naturality}).

We see that independence inside the category of commutative quantum probability spaces is fixed uniquely by natural properties which nobody wants to miss. However, comming back to the intrinsic interpretation, although the statement ``$X_2$ does not influence $X_1$'' is not symmetric under exchange of $X_1$ and $X_2$, classical independence is symmetric. In other words, within a single classical probability space there is no possibility to find random variables $X_1,X_2$ such that $X_2$ is independent of $X_1$ but not conversely. Is this realistic? In other words, are all situations in which we want to derive probabilistic results of the type that, if measurement of $X_2$ does not influence knowledge on $X_1$, then measurement of $X_1$ does not influence knowledge on $X_2$?

We consider an example. Suppose $X_1$ is an observable at instant $t_1$ and$X_2$ is an observable at instant $t_2>t_1$. We naturally expect that $X_1$ is independent of $X_2$. (A measurement at time $t_2$ cannot influence earlier times.) $X_2$ may depend on whether we measured $X_1$ or not. In the extreme case, measuring $X_1$ ``destroys'' what we know about $X_2$. (Conditioning $X_2$ on $X_1$, we obtain the constant function $\E(X_2)$ of $X_1$.)

Can we model this situation by a \hl{classical} model (i.e.\ commutative algebras)? \hl{No}, certainly not, because classical independence is symmetric. Can we model the situation by a \hl{noncommutative} model ($X_i$ are represented as operators on a Hilbert space, say)? \hl{Yes}, as follows. As before, we represent $X_2$ by the operator
$$
\id\otimes X_2~~~\text{on}~~~L^2(\R,\mu_1)\otimes L^2(\R,\mu_2).
$$
For $X_1$, however, we choose
$$
X_1\otimes 11^*~~~~~~~~~~~~~(11^*=|1\rangle\langle 1|\colon f\mapsto 1\AB{1,f}).
$$

\noindent
We find\vspace{-2ex}
\bmu{\label{monodef}
\vp_{12}\bfam{g_0(X_2)f_1(X_1)g_1(X_2)\ldots f_n(X_1)g_n(X_2)}
\\
~=~
\E\SB{g_0(X_2)}\E\SB{g_1(X_2)}\ldots\E\SB{g_n(X_2)}~\cdot~\E\bSB{(f_1\ldots f_n)(X_1)}.
}\emu
In other words, the (quantum) random variables $X_1$ and $X_2$ (more precisely, the algebras $\cA_1\otimes 11^*$ and $\id\otimes\cA_2$ generated by them) are \hl{monotone independent} in the sense of Lu \cite{Lu97} and Muraki \cite{MurN97}.

A possible interpretation: Given $f(X_1)$, measuring $g(X_2)$, the information on $X_1$ is fully preserved. Indeed, as the expression
$$
f'(X_1)g(X_2)f(X_1)
~=~
(f'f)(X_1)~\E\SB{g(X_2)}
$$
shows, a second measurment $f'(X_1)$ ``sees'' the full function $f(X_1)$ so that the product $(f'f)(X_1)$ occurs. On the other hand, the measurment of $f'(X_1)$ has completely ``nivelated'' any information about $X_2$ present after the measurment of $g(X_2)$. Provocantly, for the conditioning on $X_1$ under the condition $X_2\in B$ we might write
$$
\E(X_1|X_2\in B,X_1)
~=~
X_1
~~~~~~~~~~~
\E(X_2|X_2\in B,X_1)
~=~
1.
$$
We see that a very natural model involving only classical random variables $X_1,X_2$ leads naturally to a (noncommutative) quantum probability space.

The commutative nature of $\cA_i$ does not play a role. Allowing for noncommutative algebras we obtain immediately (quantum) \hl{tensor} independence and (quantum) \hl{monotone} independence. (Observe that the representation space for the joint distribution is the same in both cases. Only the embedding of $\cA_1$ --- unital in the tensor case, and non-unital in the monotone case --- differs.)

There are many other (quantum) independences, in the first place, Voiculescu's \cite{Voi87} famous \hl{free} independence and \hl{boolean} independence (von Waldenfels \cite{Wal73}). Additionaly, many authors consider notions of convolutions which lead to (quantum) central limit theorems as independences. However, by the following theorem we see that only the four mentioned independences fulfill Schürmann's axioms. ($q$--Convolution, for instance, is not \it{natural}; see van Leeuven and Maassen \cite{LeMa95p}.)

\stru{Theorem (Speicher \cite{Spe97}, Ben Ghorbal and Schürmann \cite{BGSc99p}, Muraki \cite{MurN02}).~}\label{1stthm}
\it{There are only two \hl{unital} and \hl{symmetric} independences, namely, \hl{tensor} and \hl{free} independence, fulfilling Schürmann's axioms. There is only one (family of) \hl{non-unital}, \hl{symmetric} indepence(s), namely \hl{boolean} independence. There is only one \hl{non-unital}, \hl{non-symmetric} independence, namely, \hl{monotone} independence.}

\stru{Remark.~}
In our example we decided that measuring $X_1$ nivelates information about $X_2$ while measuring $X_2$ leaves unaffected information about $X_1$. Of course, there are many thinkable ways how $X_2$ can influence $X_1$. However, the preceding theorem tells us that only the nivelation as described in our example can be intepreted as (quantum) independence.

\section{Conditional independence}\label{condindsec}

Let us consider a (classical) game: Two players 1,2 throw three coins, one ``fair'' coin ($p_{head}=p_{tail}$), one ``good'' coin (favouring head) and one ``bad'' (favouring tail), in the following way. First, they throw the fair coin. Depending on the result, player 1 throws the ``good'' coin and player 2 throws the ``bad'' one, or conversely. Denote by

\lf\noindent
\begin{tabular}{rl}
~~~~~~~~~~~~~$Y$&the outcome of the first throw (of the ``fair'' coin),
\\
$X_i$&the outcome of the second throw for player $i$.
\end{tabular}

\lf\noindent
Of course, $X_1$ and $X_2$ are not independent. They are, however, independent given that $Y$ has a certain value. In other words, they are \hl{conditionally} \hl{independent} (over $Y$), i.e.\ \vspace{-1ex}
\beq{ \label{cefac}
\E\bfam{f(X_1)g(X_2)\,\big|\,Y}
~=~
\E\bfam{f(X_1)\,\big|\,Y}~\E\bfam{g(X_2)\,\big|\,Y}.\vspace{-1ex}
}\eeq
Equivalently, the conditional probabilities factorize
$$
P(y,dx_1\times dx_2)
~=~
P(y,dx_1)~P(y,dx_2).
$$
We make two observations:
\begin{itemize}
\item
Given the probability spaces $(\R\times\R,\mu_{0i})$ describing each pair $(Y,X_i)$ (or, equivalently, the correponding $L^2$--spaces) it seems difficult to construct (by an abstract construction involving only the abstract strucure of probability spaces) the joint probability space $(\R\times\R\times\R,\mu)$ for $(Y,X_1,X_2)$ directly from the former.
\item
The above factorizations do not depend on the actual distribution $\nu$ of $Y$, but only on its $\sigma$--algebra, i.e.\ on the algebra $L^\infty(\R,\nu)$.
\end{itemize}
Quantum probability can help solving the problem in the first issue taking into account the second. What do we have? We have two pairs $(\cA_i,\Phi_i)$
\baln{
\cA_i
~&=~
L^\infty(\R\times\R,\mu_{0i})
&&\text{an algebra}
\\
\Phi_i\colon\cA_i\rightarrow\cA_0&=L^\infty(\R,\nu)\subset \cA_i
&&\text{a conditional expectation}
}\ealn
with $\Phi_i\bfam{f(Y)g(X_i)}=\E\bfam{f(Y)g(X_i)\big|Y}=\E\bfam{g(X_i)\big|Y}f(Y)$. Our task is to construct a joint (quantum) probability space with a conditional expectation $\Phi$ factorizing according to \eqref{cefac}. The problem: Be it the product of classical probability spaces or be it the tensor product of quantum probability spaces, in both products we end up with two copies of $Y$ (or better of the subalgebra $\cA_0$) contained in each factor.

How can we identify suitably members of the two different copies? The answer follows from the observation that by the bimodule property of conditional expectations, the factorization \eqref{cefac} does not ``see'' from which factor $\cA_i$ an element $h(Y)$ comes:
\baln{
\E\bfam{f(X_1)h(Y)g(X_2)\,\big|\,Y}
&~=~
\E\bfam{f(X_1)h(Y)\,\big|\,Y}~\E\bfam{g(X_2)\,\big|\,Y}
\\
&~=~
\E\bfam{f(X_1)\,\big|\,Y}~h(Y)~\E\bfam{g(X_2)\,\big|\,Y}
\\
&~=~
\E\bfam{f(X_1)\,\big|\,Y}~\E\bfam{h(Y)g(X_2)\,\big|\,Y}.\vspace{-1ex}
}\ealn
Therefore, we can amalgamate over $\cA_0$ and consider the \nbd{\cA_0}\nbd{\cA_0}module tensor product of the \nbd{\cA_0}\nbd{\cA_0}modules $\cA_i$, instead of the usual tensor product. (As mentioned this construction depends only on the algebra $\cA_0$ and the conditional expectations $\Phi_i$, but not on $\nu$.) Then the module tensor product
$$
\cA_1\odot\cA_2
~~:=~~
\cA_1\otimes\cA_2\text{\raisebox{-.5ex}{\Large$/$}}\text{\raisebox{-.75ex}{$a_1a_0\otimes a_2-a_1\otimes a_0a_2$}}
$$
is an algebra (with multiplication $(a_1\odot a_2)(a'_1\odot a'_2)=a_1a'_1\odot a_2a'_2$) isomorphic to (a dense subalgebra of) $L^\infty(\R\times\R\times\R,\mu)$ and $\Phi(a_1\odot a_2)=\Phi_1(a_1)\Phi_2(a_2)$ defines a conditional expectation such that
$$
\Phi\bfam{f(X_1)g(X_2)}
~=~
\E\bfam{f(X_1)\,\big|\,Y}~\E\bfam{g(X_2)\,\big|\,Y}.
$$
In other words, $\cA_1$ and $\cA_2$ are tensor independent \hl{with amalgamation} (over $\cA_0$) in the sense of Skeide \cite{Ske03a,Ske99a}. The terminology ``with amalgamation'' follows \cite{Voi95,Spe98}. Today we would prefer to say $\cA_1$ and $\cA_2$ are \hl{conditionally tensor independent}.

On the level of representation spaces the product of probability spaces was reflected by the tensor product. We now explain that for conditional independence the ``good'' representation spaces are \hl{Hilbert \nbd{\cA_0}\nbd{\cA_0}modules} and the ``correct'' product construction, once again, is the tensor product over $\cA_0$. Let
\baln{
&E_i~=~\ol{\cA_i/\sN_i}
\\
\text{with inner product}\hspace{4ex}
&\AB{x,x'}~=~\Phi_i(x^*x')
\\
\text{so that}\hspace{4ex}
&\Phi_i(a)~=~\AB{1,a1}
}\ealn
where $\sN_i$ indicates the subbimodule of length-zero elements. (This is Paschke's \cite{Pas73} GNS-construction for CP-maps applied to the simpler case of conditional expectations as in Rieffel \cite{Rie74}.) The \hl{tensor product} over $\cA_0$
\bal{\label{tpdef}
&E_1\odot E_2~=~\ol{E_1\otimes E_2/\sN}\notag
\\
\text{with inner product}\hspace{4ex}
&
\AB{x_1\odot x_2,x'_1\odot x'_2}~=~\bAB{x_2,\AB{x_1,x'_1}x'_2}
}\eal
carries a natural representation of $\cA_1\odot\cA_2$ (inherited from $\cA_1\odot\cA_2$) and we recover $\Phi$ as $\Phi(a)=\AB{1,a1}$ (where $1=1\odot 1$).

\stru{Remark.~}
Of course, for our special case in \eqref{tpdef} we might write also $\AB{x_1\odot x_2,x'_1\odot x'_2}=\AB{x_2,x'_2}\AB{x_1,x'_1}$ for the inner product, because everything commutes and $\Phi_2$ is a conditional expectation. However, if we write the definition of the inner product of a tensor product as in \eqref{tpdef}, then it works for arbitrary Hilbert bimodules. It seems to be a general feature that writing things in the ``correct'' order helps both generalizing to the noncommutative case and also understanding the classical case better.

\lf
What about conditional tensor independence for arbitrary quantum probability spaces? Here we run into serious problems. Given two algebras $\cA_i$ with a subalgebra $\cA_0$, there is only in exceptional cases a reasonable multiplication on $\cA_1\odot\cA_2$. (For instance, if $\cA_0=\sB(H)$.) Likewise, in general,
\baln{
E_1\odot E_2~&\ncong~E_2\odot E_1
&
\exists~\cA_1&\longrightarrow\cA_1\odot\id
&
\nexists~\cA_2&\longrightarrow\id\odot\cA_2.
}\ealn
This shows that, in general, conditional tensor independence does not have an analogue in quantum probability.

\stru{Remark.~}
The non-isomorphism $E_1\odot E_2\ncong E_2\odot E_1$ is at the heart of the non-possibility to construct the amalgamated analogue of tensor independence. For instance, if one of the tensor products is the null space, then certain moments, computed according to the rules in the proof of Theorem \ref{1stthm}, must be $=0$ and $\ne0$ simultaneously. At least, in the case, of von Neumann algebras plus the condition $E_1\odot E_2\cong E_2\odot E_1$, we think that it is possible to overcome the diffculties and to define the analoge of tensor independence. Notice, however, that there is no flip $x_1\odot x_2\mapsto x_2\odot x_1$ available. The embedding $\cA_2\rightarrow\sB^(E_1\odot E_2)$ is not that comming from the isomorphism $E_1\odot E_2\rightarrow E_2\odot E_1$ but it is considerably more subtle.

\lf
However, there exists always an embedding $\cA_2\rightarrow$ ``$11^*\odot\cA_2$'', where ``$11^*\odot\cA_2$'' is a not completely exact abbreviation for operators of the form
\beq{\label{projdef}
(11^*\odot a_2)(x_1\odot x_2)~:=~(1\odot\id)a_2(1^*\odot\id)(x_1\odot x_2)=~1\odot a_2\AB{1,x_1}x_2~.
}\eeq
(It is clear that in the scalar case $\cA_0=\C$  this reduces to the well-defined operator $11^*\otimes a_2$ from Section \ref{indsec}.)

\noindent
Since $(11^*\odot a_2)(a_1\odot\id)(11^*\odot a'_2)=11^*\odot a_2\Phi(a_1)a'_2$, we find
\beq{ \label{ammondef}
\Phi\bfam{a_1^{(0)}a_2^{(1)}a_1^{(1)}\ldots a_2^{(n)}a_1^{(n)}}
~=~
\Phi(a_1^{(0)})~\Phi\Bfam{a_2^{(1)}\,\Phi(a_1^{(1)})\ldots a_2^{(n)}}~\Phi(a_1^{(n)}).
}\eeq
By analogy with \eqref{monodef} we propose a new definition such that the algebras $\cA_1\odot\id$ and $11^*\odot\cA_2$ become conditionally monotone independent over $\cA_0$.

\stru{Definition.~}
Let $\cA$ be a \nbd{*}algebra with a conditional expectation $\Phi$ onto the $*$--sub\-al\-ge\-bra $\cA_0$. An ordered pair $(\cA_1,\cA_2)$ of \nbd{*}subalgebras $\cA_1,\cA_2$ containing $\cA_0$ is \hl{conditionally monotone independent} (over $\cA_0$ in $\Phi$), if they fulfill \eqref{ammondef}.

\stru{Remark.~}
Comparing with our scalar example of monotone independence, we see that $\cA_1$ and $\cA_2$ have changed their roles. The non-unitally embedded algebra (standing for events in the past of the other) is now $\cA_2$. Contrary to all good habits, the module picture shows:
$$
\text{The future is on the left of the past!}
$$

\lf
We close the discussion of conditional independence with the analogue of Theorem \ref{1stthm}. Every independence except tensor independence has an anologue conditional independence. Therfore, in Theorem \ref{1stthm} we have just to cancel tensor independence (of course, there cannot be more conditional independences, because this is not possible already in the scalar case) and arrive at the pleasant situation where every type of independence has now exactly one realization.

\stru{Theorem.~}
\it{There is exactly one \hl{unital} conditional independence, namely, conditional \hl{free} independence (which is \hl{symmetric}, automatically). There is exactly one \hl{non-unital}, \hl{symmetric} conditional independence, namely conditional \hl{boolean} independence. There is exactly one \hl{non-symmetric} conditional independence, namely, conditional \hl{monotone} independence (which is \hl{non-unital}, automatically).}

\lf\noindent
Conditional free independence was introduced by Voiculescu \cite{Voi95} and studied intensively by Speicher \cite{Spe98}. Conditional boolean independence is mentioned in Skeide \cite{Ske00} (without the simple proof of existence).

\section{Dilations}\label{dilsec}

Reversible and irreversible (quantum) dynamical systems are related by dilations. A dilation of a (small) irreversible system is a (big) reversible system such that the irreversible dynamics of the small system can be understood as a projection (conditional expecation) from the reversible dynamics of the big system onto the small subsystem. We illustrate this.
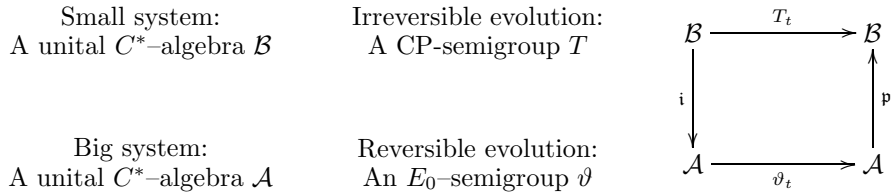
\beqn{
\parbox{10cm}{
\xymatrix{
\txt{Small system:\\{A unital $C^*$--algebra} $\cB$}
&
\txt{Irreversible evolution:\\A {CP-semigroup} $T$}
&
\cB	\ar[rr]^{T_t}	\ar[d]_\ei	&&	\cB
\\
\txt{Big system:\\{A unital $C^*$--algebra} $\cA$}
&
\txt{Reversible evolution:\\An {$E_0$--semigroup} $\vt$}
&
\cA	\ar[rr]_{\vt_t}			&&	\cA	\ar[u]_\ep
}
}
}\eeqn
A CP-semigroup is a semigroup of completely positive mappings, while an $E_0$--semi\-group is a semigroup of unital endomorphisms.

Both the upper and the lower half of the diagram are closely related to and, actually, connected by tensor product systems of Hilbert modules as introduced in Bhat and Skeide \cite{BhSk00}. We explain this in the case of a classical Markov process (joint work with L.\ Accardi to be published elsewhere).

We consider

\lf\noindent
{~}
\hfill
\begin{tabular}{ll}
$\cB=L^\infty(\R,\nu)$
&
$Y\sim\nu$,
\\
$\SB{T_t(f)}(y)=\int P_t(y,\,dx\,) f(x)$
&
$\bfam{P_t}$ a semigroup of absolutely
\\
&continuous transition functions,
\\
$\cA=L^\infty\bfam{\prod_{t\in\R_+}\R,\mu}$
&
$\mu$ the measure from Kolmogorov-Daniell
\\
&construction,
\\
$\vt_t(X_s)=X_{s+t}$
&
$\bfam{X_t}$ the associated Markov process,
\\
$\ei(Y)=X_0$
&
the canonical identification of $\cB$ with
\\
&
the subalgebra generated by $X_0$,
\\
$\ep=\E(\bullet|Y)$.
&
\end{tabular}
\hfill
{~}

\lf\noindent
Once again, what we do in the sequel does not depend on the measure $\nu$ but only on its measure type. The technical condition of $P_t(x,\bullet)$ being absolutely continuous with respect to $\nu$ assures that $T_t$ is well-defined. In practise, one will even choose $\nu$ according to this requirement, even if the distribution of $X_0$ will be different. (For instance, it could be $\delta_0$; see Remark \ref{Lrem}.)

We define a filtration. For $I\subset\R_+$ set $\cA_I=L^\infty\bfam{\prod_{t\in I}\R,\mu_I}$ and $\cA_t=\cA_{\SB{0,t}}$. Set
\baln{
E~&=~\ol{\cA/\sN}
&
E_t~&=~\ol{\cA_t/\sN_t}
&
\AB{x,x'}~&=~\ep(x^*x').
}\ealn
Then $E$ is a Hilbert \nbd{\cB}module, $E_t$ are Hilbert \nbd{\cB}\nbd{\cB}modules (left action via $X_t$, i.e.\ $f\in L^\infty(\R,\nu)$ acts as multiplication by the function $f(X_t)$) such that
\baln{
E\odot E_t~&=~E
~&~
E_s\odot E_t~&=~E_{s+t}
~&~
\vt_t(a)~&=~a\odot\id_{E_t}
\\
1\odot 1_t~&=~1
~&~
1_s\odot 1_t~&=~1_{s+t}
~&~
T_t(b)~&=~\AB{1,\vt_t\circ\ei(b)1}~=~\AB{1_t,b1_t}.
}\ealn
For instance, looking at the definition of the tensor product in \eqref{tpdef}, one checks that
\bmun{
\bfam{f_n(X_{p_n})\ldots f_1(X_{p_1})}\odot \bfam{g_m(X_{q_m})\ldots g_1(X_{q_1})}
\\
~\longmapsto~
f_n(X_{p_n+t})\ldots f_1(X_{p_1+t})g_m(X_{q_m})\ldots g_1(X_{q_1})
}\emun
defines an isomorphism $E\odot E_t\rightarrow E$. The restrictions to $E_s\odot E_t$ maps onto $E_{s+t}\subset E$ and also the left multiplication is preserved.

In particular, $E^\odot=\bfam{E_t}$ is a \hl{product system} ($E_s\odot E_t=E_{s+t}, E_0=\cB$) and $1^\odot=\bfam{1_t}$ is a \hl{unit} ($1_s\odot 1_t=1_{s+t},1_0=\U$) for $E^\odot$.

Of course, the \nbd{E_0}semigroup $\vt$ extends to all of $\sB^a(E)$ (still sending $a$ to $a\odot\id_{E_t}$). In Skeide \cite{Ske02} we showed that for every \nbd{E_0}semigroup on some $\sB^a(E)$ ($\vt_t$ being continuous in a certain topology called \hl{strict}) we can construct a (unique) product system $E^\odot$ such that $\vt_t$ can be recovered from a factorization $E=E\odot E_t$. (Acutally, we require existence of a unit vector $\xi\in E$, i.e.\ $\AB{\xi,\xi}=\U$. The general case can be treated for von Neumann algebras; see \cite{MSS02p}. The construction in \cite{Ske02} generalizes a construction for Hilbert spaces from Bhat \cite{Bha96}, while the construction in \cite{MSS02p} generalizes the construction of product systems of Hilbert spaces from Arveson \cite{Arv89}.)

\lf
We return to the construction of product systems from CP-semigroups on general $C^*$--al\-ge\-bras, and we want to see how monotone independence arises naturally. In \cite{BhSk00} we had to face a problem, namely, the possibility of a unital embedding $\cB\rightarrow\cA=\sB^a(E)$ is restricted to commutative $\cB$. Instead, we had to use $\ei(b)=1b1^*$ leading to a different filtration. Setting $\breve{\cA_0}=\ei(\cB)$, the filtration consists of the algebras $\breve{\cA_t}=\alg\CB{\vt_s(\breve{\cA_0})~(0\le s\le t)}$. While $\cA_t$ is commutative and, therefore, cannot coincide with $\sB^a(E_t)$, we find that (a suffcient closure of) $\breve{\cA_t}$ is isomorphic to $\sB^a(E_t)$ embedded as $\breve{\cA_t}=11^*\odot\sB^a(E_t)$. Also the closure of $\breve{\cA}_\infty$ is all of $\sB^a(E)$.

Neither the algebras $\cA_t$ and $\cA_{\SB{t,t+s}}$ are conditionally tensor independent nor the algebras $\breve{\cA_t}$ and $\breve{\cA\,}_{\!\SB{t,t+s}}:=\vt_t(\breve{\cA_s})$ are conditionally monotone independent. For that the dilation must be a \hl{white noise}, i.e.\ a dilation of the trivial semigroup $T_t=\id_\cB$, in yet other words, the conditional expectation $\ep=\AB{1,\bullet 1}$ must be \hl{invariant} for $\vt$. The same computations leading to \eqref{ammondef} show also the following.

\stru{Theorem.~}\label{cmithm}
\it{Let $(E,\vt,\xi)$ be a white noise, i.e.\ $\vt$ is a (strict) \nbd{E_0}semigroup on $\sB^a(E)$ leaving $\ep=\AB{1,\bullet 1}$ invariant. Set $\breve{\cA_t}=\xi\xi^*\odot\sB^a(E_t)$ (cf.\ \eqref{projdef}) and $\breve{\cA\,}_{\!\SB{t,t+s}}=\vt_t(\breve{\cA_s})$. Then for all $0<r<s<t$ the algebras $\breve{\cA\,}_{\!\SB{r,s}}$ and $\breve{\cA\,}_{\!\SB{s,t}}$ are conditionally monotone independent over $\breve{\cA_0}$.}

\stru{Remark.~}
Even if $\vt$ is not a white noise, then $\breve{\cA\,}_{\!\SB{r,s}}$ and $\breve{\cA\,}_{\!\SB{s,t}}$ are conditionally monotone independent over $\breve{\cA_s}$. This is just the (quantum) Markov property.

\stru{Rematk.~}
If $T_t=\id_\cB$ then $E^\odot$ constructed as in \cite{BhSk00} is trivial ($E_t=\cB$). However, most dilations are \hl{cocycle perturbations} of a (non-trivial, if $T$ is non-trivial) white noise, and cocycle perturbation \hl{does not change} the product system; see \cite{Ske02}.

If there is a (non-trivial) white noise, then by Theorem \ref{cmithm} there is always a (non-trivial) filtration of subalgebras $\cA_t$ such that the \hl{increment algebras} $\cA_{s,t}$ to disjoint intervals are conditionally {monotone independent}. This justifies, in particular, the name \it{white noise}.

\stru{Remark.~}
In the scalar case $\cB=\C$ all (unital) CP-semigroups are trivial. Therefore, if $\vt$ is an $E_0$--semigroup on some $\sB(\eH)$ (Arveson \cite{Arv89}) and $\xi\in\eH$ a unit vector such that $\vt_t(\xi\xi^*)$ is increasing, then the increment algebras $\vt_s(\id_\eH\otimes\sB(\eH_t))$ and $\vt_s(\om\om^*\otimes\sB(\eH_t))$ are tensor independent and monotone independent, respectively.

\stru{Remark.~}\label{Lrem}
Suppose $X=\bfam{X_t}$ is a stationary independent increment process. Thus, $Y=\bfam{Y_t}$ with $Y_t=X_t-X_0$ is a Lévy process. Then $Y_0\sim\delta_0$. Since $L^\infty(\R,\delta_0)=\C$, the associated product system consists of Hilbert spaces. It is not difficult to see that the $L^2$--space of $Y$ is $\eH=\G(L^2(\R_+,K))$ and that the associated product system is $\eH_t=\G(L^2(\SB{0,t},K))$ (Parthasarathy and Schmidt \cite{PaSchm72}). The product system $E^\odot$ of $X_t$ is intimately but non-trivially related to $\eH^\otimes$ (joint work with U.\ Franz and M.\ Schürmann; unpublished).

\lf
For a more detailed survey on product systems of Hilbert modules and the present status of their classification we refer the reader to Skeide \cite{Ske03b}.

{\small\renewcommand{\tt}{}
\newcommand{\Swap}[2]{#2#1}\newcommand{\Sort}[1]{}

}

\end{document}